%% file: agt-4-34.tex
\documentclass{gtart_h}

\input agtout

\lognumber{34}
\volumenumber{4}
\volumeyear{2004}
\papernumber{34}
\pagenumbers{781}{812}
\received{7 August 2004} 
\accepted{31 August 2004}
\published{23 September 2004}

\usepackage{amsmath,amssymb}

\usepackage[mathscr]{eucal}
\usepackage{amssymb, amsfonts, xypic}
\usepackage{enumerate}

\newtheorem{theorem}{Theorem}[section]
\newtheorem{lemma}[theorem]{Lemma}
\newtheorem{corollary}[theorem]{Corollary}
\newtheorem{proposition}[theorem]{Proposition}

\theoremstyle{definition}
\newtheorem{definition}[theorem]{Definition}
\newtheorem{hypothesis}[theorem]{Hypothesis}

\theoremstyle{remark}
\newtheorem{remark}[theorem]{Remark}

\numberwithin{equation}{section}

\newcommand{\Z}{\mathbb{Z}}
\newcommand{\A}{\mathbb{A}}

\newcommand{\C}{\mathcal{C}}

\DeclareMathOperator{\Hom}{Hom}
\DeclareMathOperator{\End}{End}

\DeclareMathOperator*{\colim}{colim}
\DeclareMathOperator*{\holim}{holim}
\DeclareMathOperator*{\Rlim}{Rlim}

\newcommand{\pros}{\textup{pro}}
\newcommand{\pro}{\textup{pro-}}
\newcommand{\ind}{\textup{ind-}}
\newcommand{\inds}{\textup{ind}}
\newcommand{\Map}{\textup{Map}}

\newcommand{\op}{\textup{op}}

\newcommand{\tl}{\tilde}
\newcommand{\map}{\rightarrow}
\newcommand{\ra}{\rightarrow}
\newcommand{\la}{\leftarrow}
\newcommand{\bdry}{\partial \Delta}
\newcommand{\bd}{\partial \Delta}

\newcommand{\smsh}{\wedge}

\newcommand{\dfn}{\textsl} 
\newcommand{\mdfn}[1]{\dfn{#1}}

\begin{document}

\title{Duality and Pro-Spectra}

\author{J. Daniel Christensen\\Daniel C. Isaksen}

\address{Dept of Math, University of Western 
Ontario, London, Ontario, Canada\\Department of 
Mathematics, Wayne State 
University, Detroit, MI 48202, USA}

\asciiaddress{Department of Mathematics, University of Western 
Ontario, London, Ontario, Canada\\and\\Department 
of Mathematics, Wayne State University, Detroit, MI 48202, USA}

\asciiemail{jdc@uwo.ca, isaksen@math.wayne.edu}

\gtemail{\mailto{jdc@uwo.ca}\qua {\rm and}\qua
\mailto{isaksen@math.wayne.edu}}

\begin{abstract}
Cofiltered diagrams of spectra, also called pro-spectra, have arisen
in diverse areas, and to date have been treated in an ad hoc manner.
The purpose of this paper is to systematically develop a homotopy
theory of pro-spectra and to study its relation to the usual homotopy
theory of spectra, as a foundation for future applications.
The surprising result we find is that our homotopy theory of
pro-spectra is Quillen equivalent to the opposite of the homotopy
theory of spectra.
This provides a convenient duality theory for all spectra, extending
the classical notion of Spanier-Whitehead duality which works well
only for finite spectra.
Roughly speaking, the new duality functor takes a spectrum to the
cofiltered diagram of the Spanier-Whitehead duals of its finite
subcomplexes.
In the other direction, the duality functor takes a cofiltered diagram
of spectra to the filtered colimit of the Spanier-Whitehead duals of 
the spectra in the diagram.
We prove the equivalence of homotopy theories by showing that both are
equivalent to the category of ind-spectra (filtered diagrams of spectra).
To construct our new homotopy theories, we prove a general existence 
theorem for colocalization model structures generalizing known results 
for cofibrantly generated model categories.
\end{abstract}


\maketitle

\section{Introduction}
\label{sctn:intro}

In recent years there have been many areas in which cofiltered diagrams 
of spectra have naturally arisen as a way to organize homotopical 
information.  For example, see the work of 
Cohen, Jones and Segal~\cite{CJS95} and Hurtubise~\cite{Hurt} 
on Floer homology theory, 
Ando and Morava~\cite{AnMo} on formal groups and free loop spaces,
and unpublished work of Dwyer and Rezk and of Arone on Goodwillie calculus.
Pro-spectra are also likely to be the target of an \'etale realization
functor on stable motivic homotopy theory~\cite{etale}.

A cofiltered diagram is called a pro-spectrum
(see Section~\ref{section:prelim-pro}) 
and in general it contains more information than its homotopy limit.
Thus it is necessary to develop a homotopy theory of pro-spectra,
and our goal is to do this systematically 
and to study its relation to the usual homotopy theory of spectra.
We expect that this framework will be useful for many applications.

Spanier-Whitehead
duality was one of the reasons for which the stable homotopy category
was invented~\cite{BG}~\cite{May-history}~\cite{SW}.
The idea is that there is a contravariant functor from
the stable homotopy category to itself that induces an equivalence
between the homotopy category of finite spectra and its own opposite.
This functor is defined by taking a spectrum $X$ to the
function spectrum $F(X, S^0)$, where $S^0$ is the sphere spectrum.

One important property of Spanier-Whitehead duality is that the double
dual of an infinite spectrum is not weakly equivalent to the original 
spectrum.  In some contexts, this is a useful property because it gives
a method for producing new interesting spectra.  On the other hand, it
is sometimes inconvenient that Spanier-Whitehead duality does not give
an equivalence between the whole stable homotopy category and its opposite.

The same situation arises in linear algebra over a
field $k$.  The functional dual induces an equivalence of the
category of finite dimensional $k$-vector spaces with its own opposite,
but it does not extend to an equivalence on the whole category
of $k$-vector spaces.  
One solution is to pass to the category of
pro-finite $k$-vector spaces.  In fact, the category of $k$-vector spaces
is equivalent to the opposite of the category of pro-finite $k$-vector spaces.

Because of the strong analogy between the stable homotopy category and
categories of chain complexes, it is natural to ask whether the same kind of
solution works for spectra.  The main result of this paper is that it does.
We define a homotopy theory for the category of pro-spectra and show that
its opposite is equivalent to the usual homotopy theory of spectra.
One use of this description of the opposite of stable homotopy theory
is to work in progress on a 
classification of polynomial functors (in the sense of Goodwillie calculus)
from spaces to spectra \cite{ABCI}.

The situation for pro-spectra is significantly more complicated 
than for pro-vector spaces 
because of the intricacies of homotopy categories.
The most naive approach is to consider pro-objects in the 
homotopy category of spectra.  The unstable version of this approach appears
in~\cite{AM} and~\cite{Sullivan}.  
Homotopy theorists have learned through countless
examples that considering diagrams in a homotopy category is usually
the wrong viewpoint.  Rather, it is better to consider commutative
diagrams in a geometric category and then study the homotopy theory
of these diagrams.  
Following this philosophy, we consider the category of pro-objects
in a geometric category of spectra and then equip it with a homotopy theory.

More precisely, we construct a model structure on pro-spectra in 
which the weak equivalences are detected by cohomotopy groups.
This model structure is contravariantly Quillen equivalent to a model structure
on the category of ind-spectra ({i.e.}, the category
of filtered systems of spectra) in which the weak equivalences
are detected by homotopy groups.  The model structure on ind-spectra
is in turn Quillen equivalent to the usual stable model structure
for spectra.

The cofibrant pro-spectra in our new model structure are easy to describe. 
They are the pro-spectra that are essentially levelwise cofibrant, that is,
they are levelwise cofibrant up to isomorphism.
The fibrant pro-spectra are only slightly harder to describe.  They are
the strictly fibrant pro-spectra that are essentially levelwise 
homotopy-finite, that is, the strictly fibrant pro-spectra (see
Section~\ref{subsection:strict})
that
are also levelwise weakly equivalent to a finite complex, up to isomorphism.
Dually, 
the fibrant ind-spectra are the essentially levelwise fibrant ind-spectra,
and the cofibrant ind-spectra are the strictly cofibrant ind-spectra that
are essentially levelwise homotopy-finite.  The importance
of finite complexes is no surprise since we are defining homotopy theories
that work well with respect to Spanier-Whitehead duality.

The description of fibrant pro-spectra in terms of homotopy-finite
spectra allows us to compute the total derived functor $\Rlim$ of 
the limit functor from pro-spectra to spectra.  For a constant pro-spectrum
$X$, $\Rlim X$ is the Spanier-Whitehead double dual of $X$.
See Remark~\ref{remark:derived-limit} for more details.

Our chief tool for establishing the appropriate homotopy theories of
pro-spectra and ind-spectra is a general existence theorem for 
a certain kind of colocalization of model categories
(see Theorem~\ref{thm:right-local}).  
This result is a generalization
of~\cite[Thm.~5.1.1]{PH} because cofibrantly generated model structures
satisfy our hypotheses, and our proof is very similar.
In our application, we begin with the strict model structure
on pro-spectra in which the weak equivalences are, up to isomorphism,
the levelwise weak equivalences (see Section~\ref{subsection:strict}).
Then we use mapping spaces {\em into} the spheres to determine the
local weak equivalences of pro-spectra; note that the model structure
on pro-spectra is a localization, not a colocalization, because we
use the dual of Theorem~\ref{thm:right-local}.
Dually, for ind-spectra,
we start with the strict structure and then use mapping spaces
{\em out of} the spheres to determine the colocal
weak equivalences.

In this paper, we need a model for spectra that has a 
well-behaved function spectrum defined on the geometric category,
not just on the homotopy category.  
We use symmetric spectra~\cite{HSS} for this model.  
Section~\ref{section:prelim-spectra} reviews the relevant ideas.
It is also possible to work entirely in the category
of $S$-modules~\cite{EKMM}.

In fact, we do not really need the full power of a general function
spectrum construction.  Rather, we only need to define function spectra
of the form $F(X, S^0)$.  It has been suggested to us that this is
probably possible on more naive categories of spectra such as
the one described in~\cite{BF}, but we have not checked the details.

It is important that the techniques that we develop here 
can be applied
to other situations involving pro-spectra.  For example, if 
$E$ is any generalized cohomology theory, then we can define a model
structure on the category of pro-spectra in which the weak equivalences
are detected by $E$-cohomology groups.  And the proof of our duality 
result extends to a proof that this model category is Quillen equivalent
to the opposite of the model category of $\End(E)$-module spectra.
Thus, we have essentially constructed all cohomological localizations
for pro-spectra.  This is relevant to work in progress
involving \'etale homotopy types and quadratic
forms over arbitrary fields \cite{DI}.
This kind of situation also occurs in \cite{cohomology}.

We assume that the reader is familiar with the language and basic results
of model categories.  The original reference is~\cite{DQ}, but
we conform to the notations and terminology of~\cite{PH}.
See also~\cite{Dwyer-Spalinski} or~\cite{Hovey}.

\subsection{Organization}
The paper is organized as follows.  We begin with the general
existence theorem for $K$-colocal model structures.  Next we review
the theory of pro-categories and ind-categories.  Then we apply
our colocalization result to
pro-categories and ind-categories.

The second part of the paper begins with a 
review of 
some details about symmetric spectra.  Next we construct and study the model
structures for pro-spectra and ind-spectra.  Finally, we prove that
the various model structures are Quillen equivalent.

\subsection{Acknowledgements}

The authors thank the SFB 343 at Universit\"at Bielefeld, Germany.  
They also thank Greg Arone for originally motivating the project
and Stefan Schwede for useful conversations.
The first author was supported by an NSERC Research Grant and
the second author was supported by an NSF Postdoctoral Research
Fellowship.

\section{Colocalizations of Model Structures}
\label{section:colocal}

In this section, we prove a general theorem about colocalization
of model structures.  Much of what appears here is very similar
to~\cite[Ch.~5]{PH}.  One important difference is that we work with
model structures that may not be cofibrantly generated.

Start with a right proper model category $\C$, that is, a model
category in which the pullback of a weak equivalence along a fibration 
is always a weak equivalence.
We refer to the cofibrations, weak equivalences, and
fibrations of $\C$ as \dfn{underlying} cofibrations, weak equivalences,
and fibrations.  
For convenience, assume that $\C$ is simplicial and write 
$\Map(\cdot, \cdot)$ for the simplicial mapping space.  In fact, the 
results of this section carry over to the non-simplicial setting, but
one must use the technical machinery of homotopy function
complexes~\cite[Ch.~17]{PH}.

Let $K$ be a set of objects of $\C$.
Since we shall only use the homotopical properties of the objects
in $K$, we may as well assume that each object in $K$ is underlying cofibrant.

\begin{definition} \label{defn:K-weak-equivalence}
A map $f\co X \map Y$ in $\C$ is a \mdfn{$K$-colocal weak equivalence} if
for each $A$ in $K$, the map
\[
\Map(A, \hat{X}) \map \Map(A, \hat{Y})
\]
is a weak equivalence of simplicial sets, where
$\hat{X} \map \hat{Y}$ is a fibrant replacement for $X \map Y$, 
that is, there is a commuting square
\[
\xymatrix{
X \ar[r] \ar[d] & \hat{X} \ar[d] \\
Y \ar[r]        & \hat{Y} }
\]
whose rows are underlying fibrant replacements.
\end{definition}

The idea is that we detect $K$-colocal weak equivalences by 
considering maps out of objects in $K$.

Observe that the choice of fibrant replacements for $X$ and $Y$ does
not matter; if the map
\[
\Map(A, \hat{X}) \map \Map(A, \hat{Y})
\]
is a weak equivalence for one choice of fibrant replacements, then
it is a weak equivalence for any other choice of fibrant replacements.
Also note that underlying weak equivalences are automatically
$K$-colocal weak equivalences.

\begin{definition} \label{defn:K-cofibration-fibration}
A map in $\C$ is a \mdfn{$K$-colocal fibration} 
if it is an underlying fibration.
A map in $\C$ is a \mdfn{$K$-colocal cofibration} if it has the left
lifting property with respect to all $K$-colocal acyclic fibrations.
\end{definition}

Because there is no difference between underlying fibrations
and $K$-colocal fibrations, we use the term ``fibration''
unambiguously for maps in either class.

We are defining a right Bousfield localization of the model category $\C$
in the sense of~\cite[Defn.~3.3.1]{PH}.
We add more weak equivalences, keep the fibrations
unchanged, and define the cofibrations to be what they must be.

Theorem~\ref{thm:right-local} states
that under some general hypotheses on $\C$, our
definitions are a model structure.
However, the two-out-of-three axiom and the retract axiom are
satisfied in general.  This follows from an inspection of the definitions.

The following results basically appear in~\cite[Ch.~5]{PH}
with minor obvious changes in the proofs.

\begin{lemma} \label{lem:Phil}
\mbox{}
\begin{enumerate}[\rm(a)]
\item The class of $K$-colocal acyclic cofibrations is the same as the class of
underlying acyclic cofibrations.
\item Let $A$ be any object of $K$.  For $n \geq 0$, the map
$i\co \bdry[n] \otimes A \map \Delta[n] \otimes A$ is a
$K$-colocal cofibration.
\item 
Let $p\co X \map Y$ be a fibration between fibrant objects
$X$ and $Y$.  Then $p$ is a $K$-colocal
acyclic fibration if and only if it has the right lifting property
with respect to every underlying cofibration
$\bdry[n] \otimes A \map \Delta[n] \otimes A$ 
in which $A$ belongs to $K$ and $n \geq 0$.
\end{enumerate}
\end{lemma}

\begin{proof}
For part (a), the proof of~\cite[Lem.~5.3.2]{PH} works word for word.
For parts (b) and (c), the proofs of~\cite[Prop.~5.2.5]{PH} and
\cite[Prop.~5.2.4]{PH} also work.  Although we do not have a set
of generating acyclic cofibrations at our disposal, we have avoided
this necessity by assuming that the map in part (c) is already a fibration.
\end{proof}

In order to prove the rest of the model structure axioms, we must
add hypotheses on $\C$.

\begin{hypothesis} \label{hyp}
Let $\C$ be a right proper simplicial model category, and let
$K$ be a set of cofibrant objects in $\C$.
Suppose that
there exists a regular cardinal $\kappa$ with the following properties.
First, each object of $K$ is $\kappa$-small relative to the underlying
cofibrations.  Second, if 
\[
X_{0} \map X_{1} \map \cdots \map X_{\beta} \map
\cdots
\]
is a $\kappa$-sequence of underlying cofibrations
and $p\co \colim_{\beta} X_{\beta} \map Y$ is a map such that
the composition $p_{\beta}\co X_{\beta} \map Y$ is a fibration
for each successor ordinal $\beta$,
then $p$ is also a fibration.
\end{hypothesis}

See Section~\ref{subsection:smallness} for a review of the notions
of smallness and $\kappa$-sequences.

The idea is that fibrations are closed under a certain kind of 
sufficiently long sequential colimit.
If $\C$ is cofibrantly generated, then we may choose $\kappa$ such that
the domains of each of the underlying
generating acyclic cofibrations as well as the objects of $K$ 
are $\kappa$-small relative to the underlying cofibrations.
Hence cofibrantly generated model categories always satisfy
Hypothesis~\ref{hyp}.
However, in our intended application to pro-spectra, 
$\C$ is not cofibrantly generated, but the above
hypothesis is still satisfied.

It is not usually a problem to find a single $\kappa$
for which each $A$ in $K$ is $\kappa$-small.  As long as each $A$ is
$\kappa_{A}$-small for some $\kappa_{A}$, 
we may choose $\kappa$ to be an upper bound for the ordinals $\kappa_{A}$.

The following lemma is {\em not} a factorization axiom for
the $K$-colocal model structure that we are constructing.  The problem
is that $K$-colocal acyclic fibrations are not detected by the 
right lifting property
with respect to the maps $\bdry[n] \otimes A \map \Delta[n] \otimes A$.
See~\cite[Ex.~5.2.7]{PH} for an example of this problem.
Using this lemma, the factorization we want follows 
from~\cite[Prop.~5.3.5]{PH}.

\begin{lemma} \label{lem:K-colocal-factorization}
Under Hypothesis~\ref{hyp},
every map $f\co X \map Y$ 
has a factorization 
into a 
$K$-colocal cofibration $i\co X \map W$ followed by 
a fibration $p\co W \map Y$ that
has the right lifting property with respect to all maps 
$\bdry[n] \otimes A \map \Delta[n] \otimes A$ for $A$ in $K$.
\end{lemma}

\begin{proof}
We use a variation on the small object 
argument~\cite[\S~10.5]{PH}.  
Let $J_{0}$ be the set of all squares
\[
\xymatrix{
\bdry[n] \otimes A \ar[r] \ar[d] & X \ar[d]^{f} \\
\Delta[n] \otimes A \ar[r] & Y   }
\]
for which $A$ belongs to $K$.
Define $Z_{0}$ to be the pushout
\[
\left( \coprod_{J_{0}} \Delta[n] \otimes A \right)
  \coprod_{  \coprod\limits_{J_{0}} \bdry[n] \otimes A  } X,
\]
and let $j_{0}\co  X \map Z_{0}$ and $q_{0}\co  Z_{0} \map Y$
be the obvious maps.

Now factor the map $q_{0}$ into an underlying acyclic cofibration
$i_{0}\co Z_{0} \map W_{0}$ followed by a fibration
$p_{0}\co  W_{0} \map Y$.  This finishes the first stage of the factorization.

We build the whole factorization by a transfinite induction of length
$\kappa$.  If $\beta$ is a limit ordinal, then set
$W_{\beta}$ to be $\colim_{\alpha < \beta} W_{\alpha}$ and set
$p_{\beta}$ to be $\colim_{\alpha < \beta} p_{\alpha}$.

On the other hand,
if $\beta$ is a successor ordinal, then define
$J_{\beta}$ to be the set of all squares
\[
\xymatrix{
\bdry[n] \otimes A \ar[r] \ar[d] & W_{\beta-1} \ar[d]^{p_{\beta-1}} \\
\Delta[n] \otimes A \ar[r] & Y   }
\]
for which $A$ belongs to $K$.
Define $Z_{\beta}$ to be the pushout
\[
\left( \coprod_{J_{\beta}} \Delta[n] \otimes A \right)
  \coprod_{  \coprod\limits_{J_{\beta}} \bdry[n] \otimes A  } W_{\beta-1},
\]
and let $j_{\beta}\co  X \map Z_{\beta}$ and 
$q_{\beta}\co  Z_{\beta} \map Y$ be the obvious maps.  
Now factor the map $q_{\beta}$ into an underlying acyclic cofibration
$i_{\beta}\co Z_{\beta} \map W_{\beta}$ 
followed by an underlying fibration
$p_{\beta}\co  W_{\beta} \map Y$.  
This finishes the $\beta$th stage of the factorization.

Transfinite induction yields a $\kappa$-sequence
\[
X \map W_{0} \map W_{1} \map \cdots \map W_{\beta} \map 
\cdots.
\]
Note that there are compatible maps $p_{\beta}\co  W_{\beta} \map Y$
for every $\beta$.
Let $W$ be $\colim_{\beta} W_{\beta}$, and let $p\co  W \map Y$
be $\colim_{\beta} p_{\beta}$.
By construction, the maps $p_{\beta}$ are fibrations
for every successor ordinal $\beta$.  
Hypothesis~\ref{hyp} implies that $p$ is a fibration.
Observe that $p$ has
the right lifting property with respect to the maps 
$\bdry[n] \otimes A \map \Delta[n] \otimes A$
as in the usual small object argument because
each $A$ is $\kappa$-small.

It remains to show that the map $i\co  X \map W$ is a $K$-colocal cofibration.
Since $K$-colocal cofibrations are defined by a left lifting property and
since $i$ is a transfinite composition of the maps
$W_{\beta} \map W_{\beta +1}$, it suffices to show that
each $W_{\beta} \map W_{\beta +1}$ is a $K$-colocal cofibration.  
But this map is the composition of $W_{\beta} \map Z_{\beta+1}$ and
$Z_{\beta +1} \map W_{\beta + 1}$.  The second map is an 
underlying acyclic cofibration by construction, 
so it is a $K$-colocal cofibration by
Lemma~\ref{lem:Phil}(a).  The first map is a cobase
change of a coproduct of maps of the form 
$\bdry[n] \otimes A \map \Delta[n] \otimes A$, each of
which is a $K$-colocal cofibration by Lemma~\ref{lem:Phil}(b).
Therefore, $W_{\beta} \map Z_{\beta + 1}$ is also a
$K$-colocal cofibration.
\end{proof}

\begin{theorem} \label{thm:right-local}
Under Hypothesis~\ref{hyp},
the $K$-colocal cofibrations, $K$-colocal weak equivalences, and 
fibrations
of Definitions~\ref{defn:K-weak-equivalence} and 
\ref{defn:K-cofibration-fibration} give a 
right proper simplicial
$K$-colocal model structure on $\C$.
\end{theorem}

\begin{remark} \label{rem:right-local}
The assumption that $\C$ be simplicial is not essential.  In fact,
the theorem is true when $\C$ is not simplicial (except that the
$K$-colocal model structure is also not simplicial).  The proof is similar
in spirit but uses the technical machinery of homotopy function complexes
\cite[Ch.~17]{PH}.
\end{remark}

\begin{proof}
From the basic axioms, 
only the lifting and factorization axioms require explanation.
One half of the lifting axiom follows from the definitions.  The
other half of the lifting axiom and one half of the factorization axiom
follow from the fact that $K$-colocal fibrations are the same
as underlying fibrations and that $K$-colocal acyclic cofibrations
are the same as underlying acyclic cofibrations (see Lemma~\ref{lem:Phil}(a)).
It remains only to factor maps into $K$-colocal cofibrations
followed by $K$-colocal acyclic fibrations.  These factorizations
can be constructed as in~\cite[Prop.~5.3.5]{PH} (see also
\cite[Prop.~5.3.4]{PH}) using Lemma~\ref{lem:K-colocal-factorization}.

This finishes the basic model structure axioms.  
Right properness follows immediately from~\cite[Prop.~3.4.4]{PH}.
The simplicial structure can be deduced from~\cite[\S~5.4.4]{PH}.
Tensors, cotensors, and mapping 
spaces are defined as in the underlying model structure.  
\end{proof}

When we consider pro-categories later, we will actually
apply the dual of Theorem~\ref{thm:right-local} to obtain
localizations rather than colocalizations.

\section{Preliminaries on Pro-Categories} 
\label{section:prelim-pro}

We now establish some background on pro-categories and 
ind-categories.
Whenever possible, we use dual definitions and arguments.  

\subsection{Pro-Categories and Ind-Categories}

\begin{definition} 
\label{defn:pro}
For a category $\C$, the category \mdfn{$\pro \C$} has objects all 
cofiltered
diagrams in $\C$, and 
$$\Hom_{\pro \C}(X,Y) = \lim_s \colim_t \Hom_{\C}
     (X_t, Y_s).$$
Composition is defined in the natural way.
\end{definition}

A category $I$ is \dfn{cofiltered} if the following conditions hold:
it is non-empty and small;
for every pair of objects $s$ and $t$ in $I$,
there exists an object $u$ together with maps $u \map s$ and
$u \map t$; and for every pair of morphisms $f$ and $g$ with the
same source and target, there exists a morphism $h$ such that $fh$ equals
$gh$.
Recall that a category is \dfn{small} if it has only
a set of morphisms.
A diagram is said to be \dfn{cofiltered} if its
indexing category is so.

Objects of $\pro \C$ are functors from
cofiltered categories to $\C$.  The original references for
pro-categories are~\cite{SGA} and~\cite{AM}.  

\begin{definition} 
\label{defn:ind}
For a category $\C$, the category \mdfn{$\ind \C$} has objects all 
filtered
diagrams in $\C$, and 
$$\Hom_{\ind \C}(X,Y) = \lim_t \colim_s \Hom_{\C}
     (X_t, Y_s).$$
Composition is defined in the natural way.
\end{definition}

The notion of a filtered diagram is dual to that of a cofiltered diagram.  
Namely, for every pair of 
objects $s$ and $t$ in the indexing category,
there exists an object $u$ together with maps $s \map u$ and
$t \map u$; for every pair of morphisms $f$ and $g$ in the
indexing category with the
same source and target, there exists a morphism $h$ such that $hf$ equals
$hg$;
and the category is small and non-empty.

Objects of $\ind \C$ are functors from
filtered categories to $\C$.  See~\cite{SGA} or~\cite{EH}
for more background on the definition of ind-categories.
For every category $\C$, the category $(\pro \C)^{\op}$ is isomorphic
to the category $\ind (\C^{\op})$.  Therefore, ind-categories can
be studied by considering pro-categories and appealing to duality.

The words ``pro-object'' and ``ind-object'' refer to objects of pro-categories
and ind-categories respectively.
A \dfn{constant} pro-object or ind-object is (isomorphic to) one indexed
by the category with one object and one identity map.
Let $c\co  \C \map \pro \C$ be the functor taking an object $X$ to the 
constant pro-object with value $X$.  Similarly, let 
$c\co  \C \map \ind \C$ take an object to its associated constant ind-object.
Note that these functors make $\C$ a full subcategory of both 
$\pro \C$ and $\ind \C$.

\subsection{Level Representations}

A \dfn{level representation} of a map of pro-objects
$f\co X \map Y$ is:
a cofiltered index category $I$;
pro-objects $\tl{X}$ and $\tl{Y}$ indexed by $I$; pro-isomorphisms
$X \map \tl{X}$ and $Y \map \tl{Y}$;
and a collection of
maps $f_{s}\co \tl{X}_s \map \tl{Y}_s$ for all $s$ in $I$
such that for all $t \map s$ in $I$, there is a commutative diagram
\[
\xymatrix{
\tl{X}_{t} \ar[d] \ar[r] & \tl{Y}_{t} \ar[d] \\
\tl{X}_{s} \ar[r]        & \tl{Y}_{s}       }
\]
and such that the maps $f_{s}$ represent a pro-map
$\tl{f}\co  \tl{X} \map \tl{Y}$
belonging to a commutative square
\[
\xymatrix{
X \ar[d] \ar[r]^{f} & Y \ar[d] \\
\tl{X} \ar[r]_{\tl{f}} & \tl{Y}     }
\]
in $\pro \C$.
In other words, a level representation of a pro-map $f$
is just a natural transformation 
such that the maps $f_{s}$ represent the element $f$ of
\[
\lim_{s} \colim_{t} \Hom_{\C} (X_{t}, Y_{s}) \cong
\lim_{s} \colim_{t} \Hom_{\C} (\tl{X}_{t}, \tl{Y}_{s}).
\]
A similar definition applies to level representations of
ind-maps.
Every map in $\pro \C$ (and in $\ind \C$) has a level representation 
\cite[Appendix 3.2]{AM}~\cite{Meyer}.

A pro-object ({resp.}, ind-object) $X$ satisfies a certain property 
\dfn{levelwise} if each $X_{s}$ satisfies that property.  Similarly,
a level representation $X \map Y$
satisfies a certain property \dfn{levelwise} if
each $X_s \map Y_s$ satisfies that property.
A pro-object, ind-object or map satisfies a certain property 
\dfn{essentially levelwise} if up to isomorphism it
satisfies that property levelwise.

\subsection{Smallness and Cosmallness}
\label{subsection:smallness}

In this section we provide some results about small objects in
ind-categories and cosmall objects in pro-categories.

Let $\lambda$ be any ordinal.  Then $\lambda$ is the partially ordered
set of all ordinals strictly less than $\lambda$.  
A \mdfn{$\lambda$-sequence} $X$ in a category $\C$ is a functor
from $\lambda$ to $\C$ such that the natural map 
$\colim_{\beta < \alpha} X_{\beta} \map X_{\alpha}$ is an isomorphism
for every limit ordinal $\alpha < \lambda$.
In other words, it is a diagram
\[
X_{0} \map X_{1} \map \cdots \map X_{\beta} \map \cdots 
\]
of length $\lambda$ with a kind of continuity condition at the limit ordinals.
Note that $X_{\beta}$ is defined only for 
$\beta < \lambda$, not for $\beta = \lambda$.
A \mdfn{$\lambda$-tower} is a contravariant functor from
$\lambda$ to $\C$
such that the natural map 
$X_{\alpha} \map \lim_{\beta < \alpha} X_{\beta}$ is an isomorphism
for every limit ordinal $\alpha < \lambda$.
In other words, it is a diagram
\[
\cdots \map X_{\beta} \map \cdots \map X_{1} \map X_{0}
\]
with a kind of continuity condition at the limit ordinals.

Recall that a cardinal $\lambda$ is \dfn{regular} if
the disjoint union
\[
\coprod_{\beta \in \mu} \nu_{\beta}
\]
is smaller than $\lambda$
whenever $\mu$ is a cardinal smaller than $\lambda$ and 
$\nu_{\beta}$ is a cardinal smaller
than $\lambda$ for every $\beta$ in $\mu$.
The relevance of regular cardinals
is that a transfinite sequence of maps whose length is a regular cardinal
has no cofinal subsequence of shorter length.

Recall that an object $X$ of a category $\C$ is \mdfn{$\kappa$-small} 
relative to a class of maps $C$ 
\cite[Defn.~10.4.1]{PH} if
for every regular cardinal $\lambda \geq \kappa$ and every
$\lambda$-sequence $Y$ such that each $Y_{\beta} \map Y_{\beta+1}$ belongs
to $C$, the natural map
\[
\colim_{\alpha} \Hom (X, Y_{\alpha}) \map \Hom(X, \colim_{\alpha} Y_{\alpha})
\]
is an isomorphism.
The idea is that every map from $X$ to $\colim Y$ must factor through
some $Y_{\alpha}$, and this factorization is unique up to refinement.
Dually, an object $X$ is \mdfn{$\kappa$-cosmall} relative to $C$ if
for every regular cardinal $\lambda \geq \kappa$ and every
$\lambda$-tower $Y$ such that each $Y_{\beta+1} \map Y_\beta$ belongs 
to $C$, the natural map
\[
\colim_{\alpha} \Hom (Y_{\alpha}, X) \map \Hom(\lim_{\alpha} Y_{\alpha}, X)
\]
is an isomorphism.
Again, the idea is that every map from $\lim Y$ to $X$ must factor through
some $Y_{\alpha}$ uniquely up to refinement.

Smallness is a common feature of objects of most familiar set-based
categories.  For example, every simplicial set is $\kappa$-small
relative to all maps
for some $\kappa$.  On the other hand, cosmallness is a rare property;
the only cosmall objects (relative to all maps) in the 
categories of topological spaces, simplicial sets, sets, groups,
sheaves, {etc.}\ are the terminal objects!
Nonetheless, pro-categories have many cosmall objects.

\begin{proposition} \label{proposition:pro-cosmall}
Let $\C$ be any category, and let $\omega$ be the first infinite ordinal.
Every constant pro-object of $\pro \C$ is $\omega$-cosmall relative
to all pro-maps.
More generally, 
if $X$ is a pro-object and $\kappa$ is an infinite regular cardinal
such that $X$ can be represented by a diagram 
with fewer than $\kappa$ morphisms,
then $X$ is $\kappa$-cosmall relative to all pro-maps.
\end{proposition}


\begin{proof}
Let $X$ be a cofiltered diagram with 
fewer than $\kappa$ morphisms, and
let $Y$ be a $\lambda$-tower in $\pro \C$ for some regular
cardinal $\lambda \geq \kappa$.
By~\cite[Thm.~7.2]{prolimits},
$\Hom_{\pros} (\lim_{\beta} Y^{\beta}, cX_{t})$ is isomorphic to
$\colim_{\beta} \Hom_{\pros} (Y^{\beta}, cX_{t})$ for every $t$.  After 
taking limits with respect to $t$, we conclude that
$\Hom_{\pros} (\lim_{\beta} Y^{\beta}, X)$ is isomorphic to
$\lim_{t} \colim_{\beta} \Hom_{\pros} (Y^{\beta}, cX_{t})$.
The indexing category for the limit is smaller than $\kappa$, and
the indexing category for the sequential colimit is at least as large as
$\kappa$.  The following lemma tells us that
the limit and colimit commute, and therefore
$\Hom_{\pros} (\lim_{\beta} Y^{\beta}, X)$ is isomorphic to 
$\colim_{\beta} \Hom_{\pros} (Y^{\beta}, X)$.
This shows that $X$ is $\kappa$-cosmall.
\end{proof}

\begin{lemma}
\label{lem:colim-lim}
Let $I$ be a cofiltering category and $\lambda$ be a regular cardinal
such that $I$ has fewer than $\lambda$ morphisms.
Suppose that $(t,\beta) \mapsto F^t_\beta$ is a set-valued functor on 
$I \times \lambda$.  Then the natural map
\[
\colim_{\beta \in \lambda} \lim_{t \in I} F^t_\beta \map
\lim_{t \in I} \colim_{\beta \in \lambda} F^t_\beta
\]
is an isomorphism.
\end{lemma}

\begin{proof}
For surjectivity, a typical element of the target is a compatible family
$\{ \overline{f^t} \}$,
where each $\overline{f^t}$ belongs to $\colim_{\beta} F^t_\beta$.
Now $\overline{f^t}$ is represented by an element $f^t_{\beta(t)}$.
Since $I$ is smaller than $\lambda$ and $\lambda$ is regular, 
there exists $\alpha < \lambda$ such that $\beta(t) \leq \alpha$
for all $t$.  In other words, each $\overline{f^t}$ is represented
by an element $f^t_{\alpha}$ of $F^t_\alpha$.  These elements do
not necessarily give us an element of $\lim_t F^t_\alpha$ because
they might not be compatible.  However, the original compatibility of the 
family $\{ \overline{f^t} \}$ tells us that 
for each map $i\co  t \map s$ of $I$, 
$f^t_\alpha$ and $f^s_\alpha$ become compatible after further refinement
to some $\alpha(i)$.  Again because $I$ is smaller than $\lambda$,
there exists a $\gamma$ such that $\gamma < \lambda$ and 
$\alpha(i) \leq \gamma$ for all maps $i\co  t \map s$ of $I$.  This shows
that the map is surjective.

For injectivity, now suppose given two elements $f$ and $g$ 
of the source that agree in the target.
We may as well assume that they are both represented by elements
of $\lim_t F^t_\beta$ for some $\beta$.  In other words, we have
two compatible families $\{ f^t_\beta \}$ and $\{ g^t_\beta \}$,
where each $f^t_\beta$ and $g^t_\beta$ belongs to $F^t_\beta$.
Since $f$ and $g$ agree in the target, we know that
$f^t_\beta$ and $g^t_\beta$ represent the same element of
$\colim_\beta F^t_\beta$ for each $t$.  
This means that $f^t_\beta$ and $g^t_\beta$
become equal after refinement to some $\beta(t)$.  As in the 
previous paragraph, there exists $\alpha < \lambda$ such that
$\beta(t) \leq \alpha$ for all $t$.  This shows that
the representatives of $f$ and $g$ are in fact equal after passing to
$\lim_t F^t_\alpha$, which shows that $f$ equals $g$.
\end{proof}

\begin{corollary} \label{corollary:pro-cosmall}
Every object of every pro-category is $\kappa$-cosmall relative
to all pro-maps for some $\kappa$.
\end{corollary}

\begin{proof}
This follows immediately from Proposition~\ref{proposition:pro-cosmall}.
\end{proof}

\begin{corollary} \label{corollary:ind-small}
Let $\C$ be any category, and let $\omega$ be the first infinite ordinal.
Every constant ind-object of $\ind \C$ is $\omega$-small relative
to all ind-maps.
More generally, 
if $X$ is an ind-object and $\kappa$ is an infinite regular cardinal
such that $X$ can be represented by a diagram 
with fewer than $\kappa$ morphisms,
then $X$ is $\kappa$-small relative to all ind-maps.
Therefore,
every object of every ind-category is $\kappa$-small relative 
to all ind-maps for some $\kappa$.
\end{corollary}

\begin{proof}
Using that the category $\ind \C$ is isomorphic to the opposite of the category
$\pro (\C^{\op})$, the corollary is dual to
Proposition~\ref{proposition:pro-cosmall} and 
Corollary~\ref{corollary:pro-cosmall}.
\end{proof}

While these results seem innocuous, they have allowed us to provide
a slick proof of Theorem~\ref{theorem:pi^*-model-structure}, 
whose proof originally occupied
many pages.

\subsection{Strict Model Structures}
\label{subsection:strict}

Let $\C$ be a proper simplicial model category, such
as the category of topological spaces, simplicial sets,
or any of the standard models for spectra.
Then the category $\pro \C$ has a 
\dfn{strict model structure} 
as originally developed in
\cite{EH}.  See~\cite{strict} for details about this.
The role of properness is to replace
the awkward niceness hypothesis of~\cite[\S 2.3]{EH}.
The \dfn{strict weak equivalences} ({resp.}, \dfn{strict cofibrations})
of $\pro \C$ are the essentially
levelwise weak equivalences ({resp.}, cofibrations).
The \dfn{strict fibrations} of
$\pro \C$ are defined by the right lifting property with respect
to strict acyclic cofibrations.
In fact, a more explicit description 
of the fibrations (up to retract) is possible
\cite[Defn.~4.2]{strict}, but we do not need this.

Dually, $\ind \C$ also has a strict structure.
The 
\dfn{strict weak equivalences} ({resp.}, \dfn{strict fibrations})
of $\ind \C$ are the
essentially levelwise weak equivalences ({resp.}, fibrations), and the 
\dfn{strict cofibrations} of
$\ind \C$ are defined by the left lifting property with respect to
acyclic fibrations.

The simplicial structure is defined in the following way.
For any two 
pro-objects $X$ and $Y$, the simplicial mapping space $\Map_{\pros}(X, Y)$ is
defined to be $\lim_{s} \colim_{t} \Map(X_{t}, Y_{s})$.
Dually, if $X$ and $Y$ are ind-objects, then the
simplicial mapping space $\Map_{\inds}(X, Y)$ is defined to be 
$\lim_{t} \colim_{s} \Map(X_{t}, Y_{s})$.

Tensors and cotensors with finite simplicial sets are defined
levelwise.  Tensors and cotensors with arbitrary simplicial sets
are defined by extending with limits and colimits from the definition
for finite simplicial sets.  See~\cite[Defn.~16.2]{prospace} for the reason
for this complication with infinite simplicial sets.

Suppose that $\C$ is stable in the sense that $\C$ is a pointed
model category in which the suspension and loops functors
form a Quillen equivalence from $\C$ to itself.
Then the strict model structures on $\pro \C$ and $\ind \C$ are again stable.
Suspensions and loops in $\pro \C$ and $\ind \C$ are defined levelwise.

We will need the following result about strict weak equivalences
and essentially levelwise properties later.

\begin{proposition}
\label{prop:ess-lev-C}
Let $\C$ be a proper simplicial model category, and let $C$ be any
class of fibrant objects of $\C$ that is closed under weak equivalences
between fibrant objects.
If $f\co  X \map Y$ is a strict weak equivalence such that $Y$ is strictly fibrant
and $X$ belongs to $C$ essentially levelwise, then $Y$
also belongs to $C$ essentially levelwise.
\end{proposition}

\begin{proof}
We may assume that $f$
is a levelwise weak equivalence and that there is an isomorphism
$X \map Z$ such that $Z$ belongs to $C$ levelwise.
We may take a level replacement for the diagram
$Z \la X \ra Y$
in such a way that $Z$ still belongs to $C$ levelwise and
$X \map Y$ is still a levelwise weak equivalence \cite[A.3.2]{AM}.
Moreover we may assume that the indexing category is a cofinite
directed set (see~\cite[\S~2.3]{strict}).
Factor the level map $X \map Z$ into a levelwise cofibration 
$X \map A$ followed by a levelwise acyclic fibration $A \map Z$.
Note that this is not a factorization in the strict model structure;
it is simply a factorization of each map $X_s \map Z_s$ in a functorial
way.

Now let $P$ be the levelwise pushout $A \coprod_X Y$, and let
$Q$ be the levelwise pushout $Z \coprod_X Y$.  
The map $A \map P$ is a cobase change of a levelwise weak equivalence
along a levelwise cofibration, so it is also a levelwise weak equivalence.

The natural map $Y \map Q$ is an isomorphism in $\pro \C$ because
it is a cobase change of an isomorphism.  Also, this map factors through
$P$.  Therefore, $Y$ is a retract of $P$.

Now take a strict acyclic cofibration $P \map P'$ such that $P'$ is
strictly fibrant.  Since $P$ is indexed by a cofinite directed set,
we may use the construction given in~\cite[Lem.~4.7]{strict} so that
$P \map P'$ is a levelwise weak equivalence and $P'$ is levelwise fibrant.  
Now we have a diagram
\[
\xymatrix{
Y \ar[r] & P \ar[d]\ar[r] & Y \\
& P' \ar@{-->}[ur]            }
\]
in which the dotted arrow exists because $P \map P'$ is a strict
acyclic cofibration and $Y$ is strictly fibrant.  Therefore,
$Y$ is also a retract of $P'$.  
The class of pro-objects that belong to $C$ essentially levelwise
is closed under retract \cite[Thm.~5.5]{prolimits}, so it suffices
to consider $P'$.

For each $s$, there is a zig-zag
\[ Z_s \la A_s \ra P_s \ra P'_s \]
of weak equivalences.  Since $Z_s$ belongs to $C$, the assumption on $C$
implies that each $P'_s$ belongs to $C$.
\end{proof}

\begin{remark}
Let $\C$ be a proper simplicial model category, and let $C$ be any
class of objects of $\C$ that is closed under weak equivalences.
Then strict weak equivalences preserve the class of objects in $\pro \C$
that belong to $C$ essentially levelwise.
The proof of this fact is basically the same as the proof of
Proposition~\ref{prop:ess-lev-C} but slightly shorter.
\end{remark}

\section{Localizations of Pro-Categories}
\label{section:pro-colocal}

In this section, we prove a
general localization result for homotopy theories of pro-categories.
The reason that we obtain localizations, not colocalizations, is that
we apply the dual of Theorem~\ref{thm:right-local}.
Let $\C$ be a proper simplicial model category.
Let $K$ be any set of
fibrant objects in $\C$, and let $cK$ be the set of constant pro-objects
$cA$ such that $A$ belongs to $K$.  

\begin{definition}\label{definition:pro-cofibration}
A map in $\pro \C$ is a \dfn{cofibration} if it is 
an essentially levelwise cofibration.
\end{definition}

These cofibrations are exactly the strict cofibrations of $\pro \C$.

\begin{definition}\label{definition:pro-weak-equivalence}
A map $f\co  X \map Y$ in $\pro \C$ is a \mdfn{$cK$-local weak equivalence}
if it induces a weak equivalence
\[
\Map_{\pros}(\tilde{X}, cA) = \colim_s \Map_{\C}(\tilde{X}_s, A) \map
\colim_t \Map_{\C}(\tilde{Y}_t, A) = \Map_{\pros}(\tilde{Y}, cA)
\]
for all $A$ in $K$, where $\tilde{X}$ and
$\tilde{Y}$ are strictly cofibrant replacements for $X$ and $Y$.
\end{definition}

\begin{definition}\label{definition:pro-fibration}
A map in $\pro \C$ is a \mdfn{$cK$-local fibration} 
if it has the right lifting property
with respect to all $cK$-local acyclic cofibrations.
\end{definition}

\begin{theorem}
\label{thm:pro-colocal}
Let $\C$ be a proper simplicial model category, and let $K$ be
any set of fibrant objects in $\C$.
Definitions~\ref{definition:pro-cofibration}, 
\ref{definition:pro-weak-equivalence}, and~\ref{definition:pro-fibration}
define a left proper simplicial model structure on $\pro \C$.
\end{theorem}

\begin{proof}
This
is an application of a dual version of Theorem~\ref{thm:right-local}
in which the objects of $K$ are required to be cosmall and cofibrations are
preserved by sequential limits.  These hypotheses are proved in
Proposition~\ref{proposition:pro-cosmall} and~\cite[Cor.~5.4]{prolimits}.
\end{proof}

We emphasize that it is {\em not} necessary that the model
category $\C$ satisfies Hypothesis~\ref{hyp} in Theorem~\ref{thm:pro-colocal}.
It is important that $\pro \C$ does satisfy this hypothesis.  This happens
by general arguments about pro-categories, not by using specific
properties of $\C$.
Theorem~\ref{thm:pro-colocal} could be stated even more generally.  There is no
need to localize with respect only to constant pro-objects, since every
pro-object is cosmall for some cardinal 
(see Corollary~\ref{corollary:pro-cosmall}).
Theorem~\ref{theorem:pi^*-model-structure} below 
is one example of the situation in Theorem~\ref{thm:pro-colocal}.
See~\cite{cohomology} for other examples.

We really do need the dual of Theorem~\ref{thm:right-local} in order
to establish the $cK$-local model structure of Theorem~\ref{thm:pro-colocal};
the dual of
\cite[Thm.~5.1.1]{PH} is not strong enough.  The problem is that 
the strict model structure for $\pro \C$ is not
fibrantly generated in general.
See~\cite[\S 5]{strict} for a proof that the strict model structure
for pro-simplicial sets is not fibrantly generated.

\begin{remark}
\label{remark:stable}
Suppose that $\C$ is a \emph{stable} model category in the sense
that the loops and suspension functors are inverse Quillen equivalences
of $\C$ with itself.  Let $K$ be a set
of fibrant objects of $\C$ such that 
for every $A$ in $K$, $\Omega A$ and $\Sigma \tilde{A}$ are weakly
equivalent to elements of $K$, where $\tilde{A}$ is a cofibrant
replacement for $A$.
In other words, $K$ is closed, up to homotopy, under suspensions and
loops.
Then it can be proved that the $cK$-local model structure is also stable.
We will not need this result, but the model structure of Theorem
\ref{theorem:pi^*-model-structure} is an example of this situation.
\end{remark}

We do not know whether the $cK$-local model structure on $\pro \C$
is always right proper, even though we are always assuming that $\C$ is 
proper.
If we assume in addition that $\C$ is stable,
then we can show that the model structure is right proper.

\begin{proposition} \label{proposition:pro-right-proper}
If $\C$ is a stable proper simplicial model category and $K$ is a set
of fibrant objects of $\C$, then
the $cK$-local model structure is right proper.
\end{proposition}

\begin{proof}
Suppose given a pullback square
\[
\xymatrix{
W \ar[r]^-{q} \ar[d]_{g} & Z \ar[d]^{f} \\
X \ar[r]_{p} & Y   }
\]
in $\pro \C$
in which $p$ is a $cK$-local fibration and 
$f$ is a $cK$-local weak equivalence.  
We want to show that 
$g$ is also a $cK$-local weak equivalence.

Let $F$ be the homotopy fibre of $p$ with respect to the
strict model structure, which is also the homotopy
fibre of $q$.  Now $\Sigma \tilde{F}$ is the homotopy cofibre 
of both $p$ and $q$, where $\tilde{F}$ is a cofibrant replacement for $F$.
We have a diagram
\[
\xymatrix{
\Map(\tilde{W}, cA) \ar[r]\ar[d] & \Map(\tilde{Z}, cA) \ar[r] \ar[d] &
  \Map(\Sigma \tilde{F}, cA) \ar[d] \\
\Map(\tilde{X}, cA) \ar[r] & \Map(\tilde{Y}, cA) \ar[r] &
  \Map(\Sigma \tilde{F}, cA)                                    }
\]
of simplicial sets in which the rows are fibre sequences.  Here
$\tilde{W}$, $\tilde{Z}$, $\tilde{X}$, and $\tilde{Y}$
are cofibrant replacements for $W$, $Z$, $X$, and $Y$,
and $A$ is any object of $K$.
The second and third vertical maps are weak equivalences,
which means that the first vertical map is also.
\end{proof}

\subsection{Fibrant Pro-Objects}
\label{subsection:pro-fibrant}

Later we shall need a more explicit description of $cK$-local fibrant
pro-objects.  This section contains this description.
We are still assuming that $\C$ is a proper simplicial model category
and that each $A$ in $K$ is fibrant.

\begin{definition}
\label{definition:K-nilpotent}
The class of \mdfn{$K$-nilpotent} objects of $\C$ is the smallest class
of fibrant objects such that:
\begin{enumerate}
\item the terminal object of $\C$ is $K$-nilpotent;
\item weak equivalences between fibrant objects preserve $K$-nilpotence;
\item and if $X$ is $K$-nilpotent, $A$ belongs to $K$, and 
$X \map A^{\bd[n]}$ is any map, then the fibre product
$X \times_{A^{\bd[n]}} A^{\Delta[n]}$ is again $K$-nilpotent.
\end{enumerate}
\end{definition}

In other words, an object is $K$-nilpotent if and only if it can
be built, up to weak equivalence,
from the terminal object by a finite sequence of base changes of
maps of the form $A^{\Delta[n]} \map A^{\bd[n]}$ with $A$ in $K$.
The terminology arises from the connection with nilpotent spaces
when $\C$ is the category of simplicial sets and $K$ is the collection
of Eilenberg-Mac Lane spaces.  See~\cite{cohomology} for details.

In this paper, the only important example occurs with $\C$ the category of
spectra and $K$ the set of spheres.  In this specific case,
we give in Proposition~\ref{prop:K-nilpotent-spectra} 
a more concrete description of $K$-nilpotent spectra.

\begin{lemma}
\label{lem:K-nilpotent}
Let $f\co  X \map Y$ be any $cK$-local weak equivalence between
cofibrant pro-objects.
Then the map 
$\Map(f,cZ)\co \Map(Y, cZ) \map \Map(X, cZ)$ is a weak equivalence
for all $K$-nilpotent objects $Z$ of $\C$.
\end{lemma}

\begin{proof}
The map $\Map(f,cZ)$ 
is a weak equivalence (even an isomorphism) when $Z = *$.
Since $X$ and $Y$ are cofibrant and $cZ$ is strictly fibrant,
the weak homotopy types of $\Map(Y, cZ)$ and $\Map(X, cZ)$ 
do not depend
on the choice of $Z$ up to weak equivalence.

It only remains to consider condition (3) of 
Definition~\ref{definition:K-nilpotent}.
Suppose for induction that $Z$ is $K$-nilpotent and that
$\Map(Y, cZ) \map \Map(X, cZ)$ 
is a weak equivalence.  Let $Z'$ be
the fibre product $Z \times_{A^{\bd[k]}} A^{\Delta[k]}$ for some object
$A$ in $K$.  Since $A^{\Delta[k]} \map A^{\bd[k]}$ is a fibration,
this fibre product is actually a homotopy fibre product, which means
that $\Map(Y, cZ')$ is the homotopy fibre product of the top row
in the diagram
\[
\xymatrix{
\Map(Y, cZ) \ar[r] \ar[d] & \Map(\bd[k] \otimes Y, cA) 
   \ar[d] & \Map(\Delta[k] \otimes Y, cA) \ar[l]\ar[d] \\
\Map(X, cZ) \ar[r] & \Map(\bd[k] \otimes X, cA) & 
    \Map(\Delta[k] \otimes X, cA), \ar[l] }
\]
and $\Map(X, cZ')$ is the homotopy fibre product of the bottom row.
The left vertical map is a weak equivalence by the induction assumption,
and the other two vertical maps are weak equivalences because the
$cK$-local model structure is simplicial and because $f$ is a
$cK$-local weak equivalence.  Thus, the induced map on homotopy
fibre products is also a weak equivalence.
\end{proof}

\begin{proposition} \label{prop:pro-fibrant}
An object $X$ of $\pro \C$ is $cK$-local fibrant 
if and only if it is strictly fibrant and essentially 
levelwise $K$-nilpotent.
\end{proposition}

That $X$ is essentially levelwise $K$-nilpotent means that
$X$ is isomorphic to a pro-object $Y$ such that
each $Y_{s}$ is $K$-nilpotent.

\begin{proof}
First suppose that $X$ is $cK$-local fibrant.
Every $cK$-local fibration is a strict fibration, so $X$ is strictly fibrant.
It remains to show that $X$ is essentially levelwise $K$-nilpotent.
By applying a functorial cofibrant replacement construction levelwise to $X$,
we get a map $\tilde{X} \map X$ which is a levelwise weak equivalence
such that $\tilde{X}$ is levelwise cofibrant (and, in particular, strictly cofibrant).
If we can show that $\tilde{X}$ is essentially levelwise $K$-nilpotent,
then we can use Proposition~\ref{prop:ess-lev-C} to conclude that
$X$ is also essentially levelwise $K$-nilpotent.
In other words, we might as well assume that $X$ is strictly cofibrant.

Consider the factorization $X \map W \map c*$ of the map
$X \map c*$ 
as described in the dual to the proof of
Lemma~\ref{lem:K-colocal-factorization},
so $W$ is $cK$-local fibrant.
Note that each $cA$ is $\omega$-cosmall by
Proposition~\ref{proposition:pro-cosmall}. 
Therefore, we may take $\kappa$
to be $\omega$, and there are no limit ordinals
in the construction of $W$.

Since $X$ is strictly cofibrant,
the dual of Lemma~\ref{lem:Phil}(c)
tells us that the map $X \map W$ is a $cK$-local acyclic cofibration.
Hence $X$ is a retract of $W$ because $X$ is $cK$-local fibrant.
The class of pro-objects having any 
property essentially levelwise is closed under retracts 
\cite[Thm.~5.5]{prolimits},
so it suffices to consider $W$.
But the class of pro-objects having any 
property essentially levelwise is also closed under cofiltered limits
\cite[Thm.~5.1]{prolimits},
so it suffices to consider each $W_{\beta}$.

Assume for induction that the pro-object $W_{\beta-1}$
is levelwise $K$-nilpotent.
We may take a level representation for the diagram
\[
\xymatrix{
& \prod_{J_{\beta}} cA^{\Delta[n]} \ar[d] \\
W_{\beta-1} \ar[r] & \prod_{J_{\beta}} cA^{\bdry[n]},     }
\]
and we construct $Z_{\beta}$ by taking the levelwise fibre product.
In fact, it is possible to construct the level representation in such 
a way that the replacement for $W_{\beta-1}$ is a diagram of objects
that already appeared in the original $W_{\beta-1}$.  This means
that the new $W_{\beta-1}$ is still levelwise $K$-nilpotent.

The construction of arbitrary products in pro-categories
\cite[Prop.~11.1]{prospace} shows that the map
\[
\prod_{J_{\beta}} cA^{\Delta[n]} 
   \map \prod_{J_{\beta}} cA^{\bdry[n]}
\]
is levelwise a finite product of maps of the form
\[
A^{\Delta[n]} \map A^{\bdry[n]}.
\]
It follows immediately that $Z_{\beta}$ is levelwise $K$-nilpotent.
Now $W_{\beta} \map Z_{\beta}$ is a levelwise weak equivalence, so 
$W_{\beta}$ is also levelwise $K$-nilpotent.
This finishes one implication.

Now suppose that $X$ is essentially levelwise $K$-nilpotent
and strictly fibrant.
We may assume that each $X_{s}$ is $K$-nilpotent.
Using the lifting property characterization of $cK$-local 
fibrant pro-objects,
we must show that the map
\[
f\co \Map(B, X) \map \Map(A, X)
\]
is an acyclic fibration of simplicial sets for every
$cK$-local acyclic cofibration $i\co  A \map B$ in $\pro \C$.
In fact, by~\cite[Prop.~13.2.1]{PH}, we may further assume that
$A$ and $B$ are cofibrant pro-objects (using that the $cK$-local
model structure is left proper).
We already know that $f$ is a fibration because of the strict
model structure.
Since $X$ is strictly fibrant, $f$ is weakly equivalent to 
$\holim_{s} f_{s}$, where $f_{s}$ is the map
\[
f_{s}\co \Map(B, cX_{s}) \map \Map(A, cX_{s}).
\]
Therefore, we need only show that each $f_{s}$ is a weak equivalence.
This is true by Lemma~\ref{lem:K-nilpotent}.
\end{proof}

\begin{remark}
\label{remark:retract}
In Definition~\ref{definition:K-nilpotent},
one might also require that the class of $K$-nilpotent objects is closed 
under retracts.  Strangely, this makes no difference in
Proposition~\ref{prop:pro-fibrant}.  The statement of that proposition is true
exactly as worded, whether or not $K$-nilpotent objects are closed
under retracts.  This surprising phenomenon arises from the surprising
way in which retracts interact with essentially levelwise properties
of pro-objects~\cite[Thm.~5.5]{prolimits}.
\end{remark}

\subsection{Ind-Categories}

All the results of this section dualize to ind-categories.  More specifically,
let $K$ be a set of cofibrant objects in a proper simplicial model
category $\C$.
A map $f\co  X \map Y$ in $\ind \C$ is a \mdfn{$cK$-colocal weak equivalence}
if it induces a weak equivalence
\[
\Map_{\inds}(cA,\hat{X}) = \colim_s \Map_{\C}(A, \hat{X}_s) \map
\colim_t \Map_{\C}(A, \hat{Y}_t) = \Map_{\inds}(cA,\hat{Y})
\]
for all $A$ in $K$, where $\hat{X}$ and
$\hat{Y}$ are strictly fibrant replacements for $X$ and $Y$.
Fibrations are just strict fibrations, and $cK$-colocal cofibrations
are defined by a lifting property.
These definitions give a right proper $cK$-colocal model structure
on $\ind \C$.  If $\C$ is stable, then this model structure is also
left proper.  An ind-object is $cK$-colocal cofibrant if and only if it is
strictly cofibrant and, up to isomorphism, it is levelwise weakly
equivalent to an object that can be built out of the initial object
by a finite sequence of cofibrant homotopy pushouts of maps of the form
$\bd[n] \otimes A \map \Delta[n] \otimes A$.

\section{Preliminaries on Spectra} 
\label{section:prelim-spectra}

Now we review some definitions and results about spectra.  
We work in
the category of symmetric spectra~\cite{HSS}.
This category has a 
proper simplicial cofibrantly generated stable model
structure.  We take this structure as the ``standard'' model structure
for spectra.  
Whenever we write ``spectrum'', we always mean ``symmetric spectrum''.
We write \mdfn{$S^{n}$} 
for a fixed cofibrant and fibrant model for the $n$th suspension 
of the sphere spectrum.

The category of symmetric spectra is closed symmetric monoidal.
This means that it has an associative commutative unital smash product 
$\smsh$ and an internal function object $F(\cdot,\cdot)$ such that
$F(Z,\cdot)$ is right adjoint to $\cdot \smsh Z$.
That is, there is a bijection between maps $X \map F(Z,Y)$ and maps
$X \smsh Z \map Y$.

We shall use the following model theoretic property of the
functor $F(\cdot, Y)$ when $Y$ is an arbitrary fixed fibrant 
spectrum~\cite[Cor.~5.3.9]{HSS}.  
Namely, $F(\cdot, Y)$ takes cofibrations to fibrations.
More precisely, if $i\co  A \map B$ is a cofibration, then
\[
F(i, Y)\co  F(B, Y) \map F(A, Y)
\]
is a fibration.

The weak equivalences in the category of spectra are defined 
in~\cite[Defn.~3.1.3]{HSS}.
The \dfn{stable homotopy category} is the category obtained by
inverting these maps, which are also called stable equivalences.
We will not repeat the definition of weak equivalence here, but using
the following definition of homotopy groups we will state an equivalent
condition below.

\begin{definition} \label{definition:homotopy-group}
For any spectrum $X$, 
let \mdfn{$\pi_{n} X$}
be the set $[ S^{n}, X]$ of maps in the stable homotopy category.
\end{definition}

When $X$ is a fibrant spectrum, $\pi_{n} X$ is isomorphic to
the traditional $n$th stable homotopy group
$\colim_{k} \pi_{n+k} X_{k}$.  More generally, we can calculate
$\pi_{n} X$ by considering the traditional $n$th stable homotopy
group of a fibrant replacement for $X$.
Weak equivalences are not defined in terms of homotopy
groups because the definition of homotopy groups depends on the
prior existence of the stable homotopy category.
Nevertheless,
the homotopy groups do detect stable weak equivalences 
of symmetric spectra in the sense that
a map $f\co  X \map Y$ is a stable equivalence if and
only if $\pi_{n} f$ is an isomorphism for every $n \in \Z$.  

Filtered colimits preserve fibrant spectra since the generating
acyclic cofibrations have compact domains~\cite[Defn.~3.4.9]{HSS}
(recall that an object $C$ of a category is compact if $\Hom(C, \cdot)$
commutes with filtered colimits).
Since filtered colimits
also preserve the traditional stable homotopy groups, it follows
that 
\[
\pi_{n} (\colim_{s} X_{s}) \cong \colim_{s} (\pi_{n} X_{s})
\]
for every filtered diagram $X$ of spectra.

\section{$\pi^{*}$-Model Structure on Pro-Spectra}
\label{section:pro-model-structure}

In this section we specialize
Theorem~\ref{thm:pro-colocal} to 
describe a model structure on the category of pro-spectra
that involves cohomotopy.
Later we compare the associated homotopy theory to ordinary
stable homotopy theory.

\begin{definition}\label{definition:pi^*-cofibration}
A map of pro-spectra is a \dfn{cofibration} if it is 
an essentially levelwise cofibration.
\end{definition}

These cofibrations are identical with strict cofibrations of pro-spectra.

Recall that the
cohomotopy group $\pi^{n}X$ of a spectrum $X$ is the group of stable weak
homotopy classes
$[X, S^{n}]$.  Thus $\pi^{*}$ is the cohomology theory represented by
the sphere spectrum $S^{0}$.  

\begin{definition}\label{definition:pi^*-weak-equivalence}
A map $f\co  X \map Y$ of pro-spectra is a \mdfn{$\pi^{*}$-weak equivalence} if it
induces an isomorphism 
$\colim_{s} \pi^{n}Y_{s} \map \colim_{t} \pi^{n}X_{t}$
for every $n \in \Z$.
\end{definition}

It is important that we are not requiring that 
$\pi^{n}X \map \pi^{n}Y$ be an ind-iso\-mor\-phism.  Non-isomorphic ind-abelian
groups may have isomorphic colimits when one allows infinitely generated
abelian groups.
It is also important that we are not using the groups 
$\pi^{n} \lim_{s} X_{s}$ or $\pi^{n} \holim_{s} X_{s}$.  In general,
very little can be said about these groups in terms of the groups
$\pi^{n} X_{s}$.  

\begin{proposition}
\label{prop:pi^*-weak-equivalence}
Let $f\co  X \map Y$ be a map of pro-spectra, 
and let $\tilde{f}\co  \tilde{X} \map \tilde{Y}$ be a strictly cofibrant 
replacement for $f$.
The following conditions are equivalent:
\begin{enumerate}
\item $f$ is a $\pi^*$-weak equivalence;
\item $\Map(\tilde{f}, cS^n)\co \Map(\tilde{Y}, cS^n) \map \Map(\tilde{X}, cS^n)$
is a weak equivalence of simplicial sets for every $n \in \Z$.
\item $\colim_t F(\tilde{Y}_t, S^0) \map \colim_s F(\tilde{X}_s, S^0)$
      is a weak equivalence of spectra.
\end{enumerate}
\end{proposition}

\begin{proof}
For any cofibrant pro-spectrum $Z$ and any $k \geq 0$,
\[
  \pi_{k} \Map(Z, cS^{n}) 
= \pi_{k} \colim_{t} \Map(Z_{t}, S^{n})
= \colim_{t} \, [Z_{t}, S^{n-k}] .
\]
The case $k=0$ tells us that condition (2) implies condition (1).

Now suppose that $f$ is a $\pi^*$-weak equivalence.
From the computation of the previous paragraph, 
we know that $\Map(\tilde{f}, cS^n)$
induces an isomorphism on all homotopy groups at the canonical
basepoints.
Since $\Map(\tilde{f}, cS^n)$ 
is weakly equivalent to the loop space $\Omega \Map(\tilde{f}, cS^{n+1})$, 
we only need to compute homotopy groups at one basepoint.
This shows that condition (1) implies condition (2).

For condition (3), note that $\pi_k \colim_t F(\tilde{X}, S^0)$
is isomorphic to $\colim_t \pi^{-k} X_t$ (and similarly for $\tilde{Y}$).
This shows that condition (3) is equivalent to condition (1).
\end{proof}

\begin{definition}\label{definition:pi^*-fibration}
A map of pro-spectra is a \mdfn{$\pi^*$-fibration} 
if it has the right lifting property
with respect to all $\pi^*$-acyclic cofibrations.
\end{definition}

\begin{theorem}\label{theorem:pi^*-model-structure}
Definitions~\ref{definition:pi^*-cofibration}, 
\ref{definition:pi^*-weak-equivalence}, and~\ref{definition:pi^*-fibration}
define a proper simplicial model structure on the category of pro-spectra.
\end{theorem}

We call this the \mdfn{$\pi^{*}$-model structure} on pro-spectra.

\begin{proof}
Proposition~\ref{prop:pi^*-weak-equivalence} tells us that we are discussing
a $cK$-localization, where $K$ is the set of spheres.  Therefore,
Theorem~\ref{thm:pro-colocal} gives us everything but 
right properness.
Since the model category of spectra is stable (see the last paragraph
of Section~\ref{section:prelim-pro}),
Proposition~\ref{proposition:pro-right-proper} gives us
right properness.
\end{proof}

Now we will identify the $\pi^*$-fibrant pro-spectra.
We say that a spectrum is \dfn{homotopy-finite} if it
is weakly equivalent to a finite complex, i.e. if its
image in the stable homotopy category
is in the thick subcategory generated by $S^{0}$.

\begin{proposition} \label{prop:K-nilpotent-spectra}
Let $K$ be the set of spheres.  A spectrum is $K$-nilpotent 
(see Definition~\ref{definition:K-nilpotent})
if and only if it is fibrant and homotopy-finite.
\end{proposition}

\begin{proof}
First suppose that $X$ is $K$-nilpotent.  We work by induction over
the number of pullbacks in the construction of $X$, noting that
the terminal object is homotopy-finite and weak equivalences preserve 
homotopy-finiteness.

To do the inductive step, assume that $X$ equals 
$Y \times_{(S^k)^{\bd[n]}} (S^k)^{\Delta[n]}$, where $Y$ is homotopy-finite.
We have to
show that $X$ is also homotopy-finite.  Each of the spectra $Y$, 
$(S^k)^{\bd[n]}$, and $(S^k)^{\Delta[n]}$ is homotopy-finite, and
the map $(S^k)^{\Delta[n]} \map (S^k)^{\bd[n]}$ is a fibration.  Therefore,
$X$ is a homotopy fibre product of three homotopy-finite spectra, so $X$
is also homotopy-finite.

Now assume that $X$ is fibrant and homotopy-finite.  
We have to show that $X$ is
$K$-nilpotent.  We induct on the number of cells in $X$.  If $X$ is
weakly contractible, then it is $K$-nilpotent by definition.
To do the inductive step, suppose that there is a fibre sequence
\[
X \map Y \map S^k,
\]
where $X$ has one more cell than $Y$ and $Y$ is $K$-nilpotent.
(This is dual to the usual way of attaching cells, but produces the 
same class of finite complexes because of Spanier-Whitehead duality.)

We claim that there is a homotopy pullback square
\[
\xymatrix{
\mbox{}* \ar[r]\ar[d] & (S^{k})^{\Delta[1]} \ar[d] \\
S^{k} \ar[r] & (S^{k})^{\bd[1]},           }
\]
where the bottom horizontal map is inclusion into the first factor.
See the following lemma for the proof.
In the diagram
\[
\xymatrix{
X \ar[r]\ar[d] & \mbox{}* \ar[r]\ar[d] & (S^{k})^{\Delta[1]} \ar[d] \\
Y \ar[r] & S^{k} \ar[r] & (S^{k})^{\bd[1]},           }
\]
the left square is also a homotopy pullback square since 
$X$ is the homotopy fibre of $Y \map S^k$.  Thus the composite 
square is also a homotopy pullback square, which means that $X$ is
$K$-nilpotent.
\end{proof}

\begin{lemma}
\label{lemma:spectra-pullback}
Let $f\co S^k \map (S^k)^{\bd[1]} = S^k \times S^k$ be the inclusion into
the first factor.  Then there is a homotopy pullback square
\[
\xymatrix{
\mbox{}* \ar[r]\ar[d] & (S^{k})^{\Delta[1]} \ar[d] \\
S^{k} \ar[r]_-f & (S^{k})^{\bd[1]}.           }
\]
\end{lemma}

\begin{proof}
Consider the pushout diagram
\[
\xymatrix{
\bdry[1]_+ \ar[r] \ar[d] & \Delta[1]_+ \ar[d] \\
\Delta[0]_+ \ar[r] & \Delta[1]    }
\]
of pointed simplicial sets, where the top horizontal arrow 
is the obvious inclusion and the left vertical arrow
takes $0$ to $0$ and takes both $1$ and the basepoint to the basepoint.
Note that $\Delta[1]$ is pointed at $1$.

If we apply the functor $\Map(\Sigma^\infty (-), S^k)$ to this diagram,
we obtain a pullback square
\[
\xymatrix{
\Map(\Sigma^\infty \Delta[1], S^k) \ar[r] \ar[d] &
  (S^k)^{\Delta[1]} \ar[d] \\
(S^k)^{\Delta[0]} \ar[r] & (S^k)^{\bdry[1]}   }
\]
of spectra in which the upper right corner is contractible.
This diagram is a homotopy pullback diagram because 
the right vertical map is a fibration.  
\end{proof}

\begin{proposition} \label{prop:pro-spectra-fibrant}
A pro-spectrum $X$ is $\pi^*$-fibrant if and only if it is essentially 
levelwise homotopy-finite and strictly fibrant.
\end{proposition}

This means that 
each $X_{s}$ has the weak homotopy type of a finite complex.
There is no compatibility requirement for the weak equivalences between
the spectra $X_{s}$ and the finite complexes.

\begin{proof}
This follows immediately from Proposition~\ref{prop:pro-fibrant}
together with Proposition~\ref{prop:K-nilpotent-spectra}.
\end{proof}

\begin{remark} \label{remark:derived-limit}
Using the above description of $\pi^*$-fibrant objects,
it is possible to describe explicitly the total derived functor $\Rlim$
of the limit functor from pro-spectra to spectra.
This functor is computed by taking the limit of a $\pi^*$-fibrant replacement.
With the strict structure on pro-spectra, $\Rlim$ is just the homotopy
limit functor.  However, with the $\pi^*$-model structure, $\Rlim$ is related
to Spanier-Whitehead duality.

Let $X$ be a spectrum; we shall calculate $\Rlim (cX)$.  Take
a $\pi^*$-fibrant replacement $\hat{X}$ for $cX$.  From
Proposition~\ref{prop:pro-fibrant}, we know that each spectrum
in $\hat{X}$ is homotopy-finite; thus, $\hat{X}$ is levelwise
weakly equivalent to $F( F(\hat{X}, S^0), S^0 )$.
Therefore, $\Rlim (cX)$ is weakly equivalent to
$\holim_s F( F(\hat{X}_s, S^0), S^0 )$, which is equivalent to
$F( \colim_s F(\hat{X}_s, S^0), S^0 )$.
A computation of homotopy groups shows that
$\pi_n \colim_s F(\hat{X}_s, S^0)$ equals $\pi^{-n} X$, so
$\colim_s F(\hat{X}_s, S^0)$ is weakly equivalent to the ordinary
Spanier-Whitehead dual $F(X, S^0)$ of $X$.  Thus, we have shown that
$\Rlim (cX)$ is equivalent to $F(F(X, S^0), S^0)$.

A similar analysis shows that if $X$ is an arbitrary pro-spectrum,
then $\Rlim X$ is equivalent to $F( \colim_s F(X_s, S^0), S^0)$,
{i.e.}, the Spanier-Whitehead dual of the colimit of the 
levelwise Spanier-Whitehead dual of $X$.
\end{remark}

\section{$\pi_{*}$-Model Structure on Ind-Spectra}

Now we proceed to ind-spectra.
All of the following definitions and results are dual to analogous
results in the previous section.  We skip the proofs
because they are no different.

\begin{definition}\label{definition:ind-fibration}
A map of ind-spectra is a \dfn{fibration} if it is an
essentially levelwise fibration.
\end{definition}

These fibrations are identical with strict fibrations of ind-spectra.

\begin{definition}\label{definition:ind-weak-equivalence}
A map of ind-spectra $X \map Y$ is a \mdfn{$\pi _{*}$-weak equivalence} if
for every $n \in \Z$, the map 
$\colim_{s} \pi _{n}X_{s} \map \colim _{t}\pi _{n}Y_{t}$
is an isomorphism.
\end{definition}

\begin{proposition}
\label{prop:ind-weak-equivalence}
Let $f\co  X \map Y$ be a map of ind-spectra, 
and let $\hat{f}\co  \hat{X} \map \hat{Y}$ be a strictly fibrant replacement
for $f$.
The following conditions are equivalent:
\begin{enumerate}
\item $f$ is a $\pi_*$-weak equivalence;
\item $\Map(cS^n, \hat{f})\co \Map(cS^n, \hat{Y}) \map \Map(cS^n, \hat{X})$
is a weak equivalence of simplicial sets for every $n \in \Z$;
\item $\colim_s X_s \map \colim_t Y_t$ is a stable weak equivalence
of spectra.
\end{enumerate}
\end{proposition}

\begin{definition}\label{definition:ind-cofibration}
A map of ind-spectra is a \dfn{$\pi_*$-cofibration} 
if it has the left lifting property
with respect to all $\pi_*$-acyclic fibrations.
\end{definition}

\begin{theorem}\label{theorem:ind-model-structure}
Definitions~\ref{definition:ind-fibration}, 
\ref{definition:ind-weak-equivalence}, and~\ref{definition:ind-cofibration}
define a proper simplicial model structure on the category of ind-spectra.
\end{theorem}

We call this the \mdfn{$\pi_{*}$-model structure} on ind-spectra.

\begin{proof}
Everything but left properness 
is an application of the dual version of Theorem~\ref{thm:pro-colocal}.
Left properness follows in a manner dual to the proof of 
Proposition~\ref{proposition:pro-right-proper}.
\end{proof}

\begin{proposition} 
Consider the smallest class of cofibrant spectra such that:
\begin{enumerate}
\item $*$ belongs to the class;
\item the class is closed under weak equivalences between cofibrant spectra;
\item and if $X$ belongs to the class and $\bd[n] \otimes S^k \map X$ 
is any map, then the pushout 
$\Delta[n] \otimes S^k \coprod_{\bd[n] \otimes S^k} X$ also belongs to
the class.  
\end{enumerate}
This class coincides with the class of cofibrant homotopy-finite spectra.
\end{proposition}

\begin{proposition} \label{prop:ind-spectra-cofibrant}
An ind-spectrum $X$ is $\pi_*$-cofibrant if and only if it is essentially 
levelwise homotopy-finite and strictly cofibrant.
\end{proposition}

\section{Comparison of Homotopy Theories}
\label{section:spectra-comparison}

This section contains the main results of this paper.  
Namely, the homotopy category of pro-spectra is the opposite of the
ordinary stable homotopy category.
First, we study the homotopy theory of ind-spectra.

\begin{theorem} \label{theorem:compare-ind-spectra}
The constant functor $c$ from spectra to ind-spectra is right adjoint to
the functor $\colim$.  These functors 
form a Quillen equivalence when considering the 
$\pi_{*}$-model structure on ind-spectra.
\end{theorem}

\begin{proof}
Let $X$ be an ind-spectrum, and let $Y$ be a spectrum.
By direct calculation, $\Hom_{\inds} (X, cY) \cong \Hom(\colim X, Y)$.
Thus $c$ and $\colim$ are adjoint.

The functor $c$ preserves fibrations and weak equivalences.  
Therefore, $c$ and $\colim$ are a Quillen pair.

By Proposition~\ref{prop:ind-weak-equivalence},
$X \map cY$ is a $\pi_{*}$-weak equivalence if and only if
$\colim_{s} X_{s} \map Y$ is a weak equivalence of spectra.
Hence, $c$ and $\colim$ form a Quillen equivalence.
\end{proof}

According to Theorem~\ref{theorem:compare-ind-spectra}, 
the homotopy category of ind-spectra is equivalent to
the ordinary stable homotopy category.  
We suspect that this model structure is not cofibrantly generated,
but we have not been able to prove it.

By Proposition~\ref{prop:ind-spectra-cofibrant},
every ind-spectrum $X$ is $\pi_{*}$-weakly equivalent to an ind-spectrum
whose objects are finite cell complexes. 
Therefore, the ordinary stable homotopy category is equivalent to
a homotopy category of ind-(finite cell complexes), but this latter homotopy
category does not arise from a model structure.

The next step is to compare the categories of pro-spectra and ind-spectra.

\begin{lemma}\label{lemma:function-spectrum-adjoint}
The contravariant functor $F(\cdot,Y)$ from spectra to spectra is its own
adjoint.
\end{lemma}

\begin{proof}
A map $X \map F(Z, Y)$ corresponds to a map
$X \smsh Z \map Y$.  This corresponds to a map $Z \map F(X, Y)$ because
$X \smsh Z$ and $Z \smsh X$ are isomorphic.
\end{proof}

Let $Y$ be a fixed spectrum.
By acting levelwise,
the functor $F(\cdot, Y)$  induces a contravariant functor from
pro-spectra to ind-spectra.  It also induces a contravariant functor
from ind-spectra to pro-spectra.

\begin{proposition} \label{proposition:function-spectrum-adjoint}
Let $Y$ be an arbitrary fixed spectrum.
The contravariant functors $F(\cdot,Y)$ from pro-spectra to ind-spectra 
and from ind-spectra to pro-spectra are adjoint in the sense
that
\[
\Hom_{\pros} ( X, F(Z, Y) ) \cong 
\Hom_{\inds} ( Z, F(X, Y) )
\]
for every pro-spectrum $X$ and every ind-spectrum $Z$.
\end{proposition}

\begin{proof}
This follows from direct computation and
Lemma~\ref{lemma:function-spectrum-adjoint}.
\end{proof}

\begin{proposition} \label{proposition:function-spectrum-Quillen-pair}
Let $Y$ be a homotopy-finite fibrant spectrum.
Then the contravariant functors $F(\cdot,Y)$ from pro-spectra to ind-spectra 
and from ind-spectra to pro-spectra are a Quillen pair between
the $\pi^{*}$-model structure on pro-spectra and the opposite of
the $\pi_{*}$-model structure on ind-spectra.
\end{proposition}


\begin{proof}
We already know that the functors are an adjoint pair by
Proposition~\ref{proposition:function-spectrum-adjoint}.
In order to show that they are a Quillen pair, we must prove
that $F(\cdot, Y)$ takes cofibrations 
({resp.}, $\pi^*$-acyclic cofibrations)
of pro-spectra to fibrations ({resp.}, $\pi_*$-acyclic fibrations)
of ind-spectra.

Since $Y$ is fibrant, $F(\cdot, Y)$ takes cofibrations of spectra to
fibrations of spectra.  Therefore, $F(\cdot, Y)$ takes levelwise 
cofibrations of pro-spectra to levelwise fibrations of ind-spectra.
It follows that $F(\cdot, Y)$ takes essentially levelwise cofibrations 
of pro-spectra to essentially levelwise 
fibrations of ind-spectra.

Now let $i\co A \map B$ be a $\pi^*$-acyclic cofibration of pro-spectra.
Fix $k \in \Z$ and
let $Z$ be a fibrant model for the spectrum $\Omega^{k} Y$.  Since
$Z$ is again homotopy-finite, the constant pro-spectrum $cZ$
is $\pi^*$-fibrant.  Therefore, the map
$\Map(B, cZ) \map \Map(A, cZ)$ is an acyclic fibration of simplicial sets.  
In particular, the map 
$\pi_0 \colim_s \Map(B_s, Z) \map \pi_0 \colim_t \Map(A_t,Z)$
is an isomorphism.
This means that 
$\pi_0 \colim_s F(B_s, Z) \map \pi_0 \colim_t F(A_t,Z)$ is also an
isomorphism (because $\pi_0 F(C,D) = \pi_0 \Map(C,D)$ for any
spectra $C$ and $D$ and because $\pi_0$ commutes with filtered colimits).

Since $\pi_0 \colim_s F(B_s,Z)$ is isomorphic to
$\pi_k \colim_s F(B_s,Y)$ (and similarly for $A$), it follows that
$\colim_s F(B_s, Y) \map \colim_t F(A_t,Y)$ is a weak equivalence
of spectra.  Proposition~\ref{prop:ind-weak-equivalence}
tells us that the map $F(i,Y)\co  F(B,Y) \map F(A,Y)$ is a
$\pi_*$-weak equivalence of ind-spectra.
\end{proof}

\begin{theorem} \label{theorem:function-spectrum-Quillen-equivalence}
The contravariant functors $F(\cdot,S^{0})$ from pro-spectra to ind-spectra 
and from ind-spectra to pro-spectra form a Quillen equivalence between
the $\pi^{*}$-model structure on pro-spectra and the opposite of
the $\pi_{*}$-model structure on ind-spectra.
\end{theorem}

\begin{proof}
We already showed in Proposition
\ref{proposition:function-spectrum-Quillen-pair} that the functors
are a Quillen pair.  
Note that the hypothesis of 
Proposition~\ref{proposition:function-spectrum-Quillen-pair}
is satisfied because $S^{0}$ is homotopy-finite.

Let $X$ be a cofibrant pro-spectrum, let $Z$ be a cofibrant ind-spectrum,
and let $f\co X \map F(Z, S^{0})$ be a map of pro-spectra.
Our goal is to show that $f$ is a
$\pi^*$-weak equivalence if and only if
the adjoint map $Z \map F(X,S^0)$ is a $\pi_*$-weak equivalence
of ind-spectra.

By Proposition~\ref{prop:pi^*-weak-equivalence},
the map $f$ is a $\pi^*$-weak equivalence if and only if the map
\[
\colim_s F(\tilde{F}(Z_s, S^0),S^0) \map \colim_t F(X,S^0)
\]
is a weak equivalence of spectra.
Here $\tilde{F}(C,D)$ refers to a cofibrant replacement for the function
spectrum $F(C,D)$.

By Proposition~\ref{prop:ind-spectra-cofibrant}, we may assume that
each $Z_{s}$ is homotopy-finite. 
Since the Spanier-Whitehead double dual of a finite complex is itself, the map
$Z \map F( \tilde{F}(Z, S^{0}), S^{0})$
is a levelwise weak equivalence.  
In particular, the map 
\[
\colim_s Z_s \map \colim_s F(\tilde{F}(Z_s,S^0),S^0)
\]
is a weak equivalence of spectra.

The previous two paragraphs imply that $f$ is a weak equivalence if and 
only if the composition
\[
\colim_s Z_s \map \colim_s F(\tilde{F}(Z_s,S^0),S^0) \map \colim_t F(X,S^0)
\]
is a weak equivalence of spectra.  
By Proposition~\ref{prop:ind-weak-equivalence}, this last map is a weak
equivalence if and only if the map $Z \map F(X,S^0)$ is a $\pi_*$-weak
equivalence of ind-spectra.
\end{proof}

\begin{corollary}[Main result]
The category of pro-spectra with its $\pi^{*}$-model structure
is Quillen equivalent to the opposite of the category of symmetric 
spectra with its usual stable model structure.
The equivalence is a composite of two Quillen pairs going in
opposite directions:
\[ \textup{pro-spectra} \map \textup{ind-spectra}^{op} \la \textup{spectra}^{op} . \]
We have indicated the directions of the left adjoints.
\end{corollary}

\Addressesr

\end{document}

%% file: agtout.tex
\def\ifplaintex{\expandafter\ifx\csname documentclass\endcsname\relax}

\def\gtp{{\mathsurround=0pt\it $\cal G\mskip-2mu$eometry \&\ 
$\cal T\!\!$opology $\cal P\!$ublications}}  

\def\Addressesr{\bigskip
{\small \parskip 0pt \leftskip 0pt \rightskip 0pt plus 1fil \def\\{\par}
\sl\theaddress\par
\medskip
\rm Email:\stdspace\tt\theemail\hfill\rm Received:\qua\receiveddate \par}}

\def\recd{{\small Received:\qua\receiveddate\ifx\reviseddate\relax
\else\qquad Revised:\qua\reviseddate\fi\par}} 

\def\lognumber#1{\def\thelognumber{#1}}
\def\volumenumber#1{\def\thevolumenumber{#1}}
\def\volumeyear#1{\def\thevolumeyear{#1}}
\def\papernumber#1{\def\thepapernumber{#1}}
\def\pagenumbers#1#2{\def\startpage{#1}\def\finishpage{#2}}
\def\published#1{\def\publishdate{#1}}

\def\received#1{\def\receiveddate{#1}}

\def\accepted#1{\def\accepteddate{#1}}

\def\asciiaddress#1{\def\theasciiaddress{#1}}
\def\asciiemail#1{\def\theasciiemail{#1}}

\let\\\par\let\thelognumber\relax\let\thevolumenumber\relax
\let\thepapernumber\relax\let\thevolumeyear\relax\let\startpage\relax
\let\finishpage\relax\let\publishdate\relax\let\receiveddate\relax
\let\reviseddate\relax\let\accepteddate\relax\let\theasciititle\relax
\let\theasciiauthors\relax\let\theasciiaddress\relax
\let\theasciiabstract\relax

\let\theasciiemail\relax

\ifplaintex
\font\logobig=cmssbx10 scaled 3836
\font\logomed=cmssbx10 scaled 2557
\else
\font\logobig=cmssbx10 scaled 4200
\font\logomed=cmssbx10 scaled 2800
\fi

\long\def\makeagttitle{   
\count0=\startpage
\agt\hfill      
\hbox to 45truept{\vbox to 0pt{\vglue -13truept{\logomed A\kern -.37em{\logobig 
T}\kern -.38em G}\vss}\hss}
\break
{\small Volume \thevolumenumber\ (\thevolumeyear)
\startpage--\finishpage\nl
Published: \publishdate}

\vglue .25truein

{\parskip=0pt\leftskip 0pt plus
1fil\def\\{\par\smallskip}{\Large\bf\thetitle}\par\medskip} \vglue
0.05truein

{\parskip=0pt\leftskip 0pt plus 1fil\def\\{\par}{\sc\theauthors}
\par\medskip}%
 
\vglue 0.03truein 

{\small\leftskip 25truept\rightskip 25truept{\bf Abstract}\stdspace\theabstract

{\bf AMS Classification}\stdspace\theprimaryclass
\ifx\thesecondaryclass\relax\else; \thesecondaryclass\fi\par
{\bf Keywords}\stdspace \thekeywords\par}\vglue 7truept

}   

\ifplaintex
\font\phead=cmsl9 scaled 950
\font\pnum=cmbx10 scaled 913
\font\pfoot=cmsl9 scaled 950
\headline{\vbox to 0pt{\vskip -4.5mm\line{\small\phead\ifnum
\count0=\startpage ISSN 1472-2739 (on-line) 1472-2747 (printed)
\hfill {\pnum\folio}\else\ifodd\count0\def\\{ }%
\ifx\theshorttitle\relax\thetitle\else\theshorttitle\fi\hfill{\pnum\folio}
\else\def\\{ and }{\pnum\folio}\hfill\ifx\theshortauthors\relax\theauthors
\else\theshortauthors\fi\fi\fi}\vss}}
\footline{\vbox to 0pt{\vglue 0mm\line{\small\pfoot\ifnum\count0=\startpage
\copyright\ \gtp\hfill\else
\agt, Volume \thevolumenumber\ (\thevolumeyear)\hfill\fi}\vss}}
\else
\font=cmsl9 scaled 1050
\font\lnum=cmbx10 
\font=cmsl9 scaled 1050
\makeatletter
\def\@oddhead{{\small\lhead\ifnum\count0=\startpage ISSN 1472-2739 
(on-line) 1472-2747 (printed)\hfill {\lnum\number\count0}\else\ifodd\count0
\def\\{ }\ifx\theshorttitle\relax \thetitle \else\theshorttitle\fi\hfill
{\lnum\number\count0}\else\def\\{ and }{\lnum\number\count0}
\hfill\ifx\theshortauthors\relax 
\theauthors\else\theshortauthors\fi\fi\fi}}\def\@evenhead{\@oddhead}
\def\@oddfoot{\small\lfoot\ifnum\count0=\startpage\copyright\ \gtp\hfill\else
\agt, Volume \thevolumenumber\ (\thevolumeyear)\hfill\fi}
\def\@evenfoot{\@oddfoot}
\makeatother
\fi
\let\maketitlepage\makeagttitle
\let\makeshorttitle\maketitlepage
\let\maketitle\maketitlepage

\newwrite\gtoutfile
\long\gdef\makeheadfile{  
{\def\\{, }\def\s{ }
\immediate\openout\gtoutfile head.xxx
\immediate\write\gtoutfile{Proxy-for: \ifx\theasciiauthors\relax
\theauthors\else\theasciiauthors\fi\s<\ifx\theasciiemail\relax\theemail\else\theasciiemail\fi>}
\immediate\write\gtoutfile{\noexpand\\}
\immediate\write\gtoutfile{Authors: \ifx\theasciiauthors\relax
\theauthors\else\theasciiauthors\fi}
{\def\\{ }\immediate\write\gtoutfile{Title: \ifx\theasciititle\relax
\thetitle\else\theasciititle\fi}}
\immediate\write\gtoutfile{Subj-class: GT or SG, GR etc}
\immediate\write\gtoutfile{MSC-class: \theprimaryclass\ifx\thesecondaryclass\relax\else, \thesecondaryclass\fi}
\immediate\write\gtoutfile{Journal-ref: Algebr. Geom. Topol. \thevolumenumber\s
(\thevolumeyear) \startpage-\finishpage}
\immediate\write\gtoutfile{Comments: Published by Algebraic and
Geometric Topology at}
\immediate\write\gtoutfile{\s\s\s  http://www.maths.warwick.ac.uk/agt/AGTVol\thevolumenumber/agt-\thevolumenumber-\thepapernumber.abs.html}
\immediate\write\gtoutfile{\noexpand\\}
\immediate\write\gtoutfile{}
\ifx\theasciiabstract\relax
\immediate\write\gtoutfile{\theabstract}\else
\immediate\write\gtoutfile{\theasciiabstract}\fi
\immediate\write\gtoutfile{}
\immediate\write\gtoutfile{\noexpand\\}
\immediate\write\gtoutfile{}
\immediate\closeout\gtoutfile}}  

\def\maketitlepage{\makeagttitle\makeheadfile}
\let\makeshorttitle\maketitlepage
\let\maketitle\maketitlepage